\documentclass[a4paper, 11pt]{amsart}
\usepackage{amsfonts, amsmath, amssymb, amsthm,mathtools, url,dsfont}
\usepackage{graphicx}
\usepackage[english] {babel}
\usepackage{enumitem}
\usepackage{array}
\usepackage{cellspace}
\usepackage{color}
\newtheorem{theorem}{Theorem}

\newtheorem{lemma}{Lemma}

\begin{document}

\begin{abstract}
We prove that there are arbitrarily large values of $t$ such that $|\zeta(1+it)| \geq e^{\gamma} (\log_2 t + 
\log_3 t) + \mathcal{O}(1)$. This essentially matches the prediction for the optimal lower bound in a conjecture of Granville and 
Soundararajan. Our proof uses a new variant of the ``long resonator'' method. While earlier implementations of this 
method crucially relied on a ``sparsification'' technique to control the mean-square of the 
resonator function, in the present paper we exploit certain self-similarity properties of a specially 
designed resonator function.
\end{abstract}

\title[]{Extreme values of the Riemann zeta function on the 1-line}

\author{Christoph Aistleitner} 
\address{Christoph Aistleitner, Institute of Analysis and Number Theory, TU Graz, Austria}
\email{aistleitner@math.tugraz.at}

\author{Kamalakshya Mahatab} 
\address{Kamalakshya Mahatab, Department of Mathematical Sciences, NTNU Trondheim, Norway}
\email{accessing.infinity@gmail.com}

\author{Marc Munsch} 
\address{Marc Munsch, Institute of Analysis and Number Theory, TU Graz, Austria}
\email{munsch@math.tugraz.at}

\thanks{The first and third author are supported by the Austrian Science Fund (FWF), project Y-901. This paper was 
written while the second author was a visitor at the TU Graz for one month; this visit was also supported by FWF 
project Y-901. Furthermore, this work was carried out during the tenure of an ERCIM ``Alain Bensoussan'' Fellowship 
of the second author.}

\maketitle

\section{Introduction}

Improving earlier work of Littlewood, in 1972 Levinson~\cite{lev} proved that
\begin{equation} \label{lev}
|\zeta(1+it)| \geq e^\gamma \log_2 t + \mathcal{O}(1)
\end{equation}
for arbitrarily large $t$.\footnote{Throughout this paper, we write $\zeta(s)$ for the Riemann zeta function and 
$\log_j$ for the $j$-th iterated logarithm.} This was further improved by Granville and Soundara\-rajan~\cite{gse}, who 
in 2006 established the lower bound
$$
\max_{t \in [0,T]} |\zeta(1+it)| \geq e^\gamma \left( \log_2 T + \log_3 T - \log_4 T + \mathcal{O} (1) \right)
$$
for arbitrarily large $t$. Their result also gives bounds for the measure of those $t \in [0,T]$ for which 
$|\zeta(1+it)|$ is of this size. Also in~\cite{gse}, Granville and Soundararajan state the conjecture that actually
\begin{equation} \label{conj}
\max_{t \in [T,2T]} |\zeta(1+it)| = e^\gamma \left( \log_2 T + \log_3 T + C_1 + o (1) \right),
\end{equation}
and even give a conjectural value of the constant $C_1$.\\

The proofs of Levinson and of Granville and Soundararajan rely on estimates for high moments of the zeta function and 
on Diophantine approximation arguments, respectively. Using a different method, the so-called \emph{resonance method}, 
Hilberdink~\cite{hilb} re-established~\eqref{lev}. This method can be traced back to a paper of Voronin 
\cite{voron}, but it was developed independently and significantly refined by Hilberdink and by Soundararajan 
\cite{sound} about 10 years ago. Roughly speaking, the functional principle 
of this method is to find a function $R(t)$ such that
$$
I_1 := \int_0^T \zeta (\sigma + it) |R(t)|^2 ~dt
$$
is ``large'', whereas
$$
I_2:=\int_0^T |R(t)|^2 ~dt
$$
is ``small''. Then the quotient $|I_1|/I_2$ is a lower bound for the maximal value of $|\zeta(\sigma+it)|$ in the 
range $t \in [0,T]$. The ``resonator'' $R(t)$ is chosen as a Dirichlet polynomial, and experience shows that it is often 
suitable to choose a function $R$ with multiplicative coefficients which can be written as a finite Euler product. This 
method can be implemented relatively easily if the length of $R$ is bounded by a small power of $T$, since then $I_2$ 
is 
the square-integral of a sum of essentially orthogonal terms. However, it is desirable to be able to control significantly longer resonator functions as well, and  
in~\cite{aist} a method was developed which allows the implementation of ``long'' resonators of length roughly $T^{\log_2 T}$. 
For $\sigma \in (1/2,1)$ this method allowed to recapture Montgomery's~\cite{MR0460255} lower bounds for extreme values 
of $|\zeta (\sigma + it)|$ by means of the resonance method, while on the critical line Bondarenko and Seip~\cite{bs1} 
used a ``long resonator'' to obtain lower bounds for extreme values of $|\zeta(1/2+it)|$ which even surpassed those established before by 
different methods.\footnote{From the work of Bondarenko and Seip it is also visible why the ``long resonator'' does not 
give any essential improvement of Montgomery's results in the case $\sigma \in (1/2,1)$, other than better values for 
the involved constants. See~\cite{bs1,bs2}.} In the present paper we will adapt the ``long resonator'' argument to the 
case $\sigma=1$, and prove the following theorem. 

\begin{theorem} \label{th1}
There is a constant $C$ such that
$$
\max_{t \in [\sqrt{T},T]} |\zeta(1+it)| \geq e^\gamma \left( \log_2 T + \log_3 T - C \right)
$$
for all sufficiently large $T$.
\end{theorem}

Note that our theorem is in accordance with the conjecture of Granville and Soundararajan in equation~\eqref{conj}. 
However, the conjecture is much stronger than the theorem. First of all, our theorem only gives a lower bound 
with error $\mathcal{O}(1)$, while the conjecture gives an asymptotic equality with error $o(1)$. Thus a proof of the 
full conjecture would also require a vast improvement of the upper bound $|\zeta(1+it)| \ll (\log t)^{2/3}$ of 
Vinogradov. Additionally, in the conjecture the range for $t$ is $[T,2T]$, while our theorem requires the larger range 
$[\sqrt{T},T]$. The requirement for such a longer range is typical for applications of the ``long resonator'', and 
appears in~\cite{aist} and~\cite{bs1} as well. The range $[\sqrt{T},T]$ in our theorem could be reduced to 
$[T^{\theta},T]$ for $\theta<1$, at the expense of replacing $C$ by some other constant $C(\theta)$.\\

Before turning to the proof of Theorem~\ref{th1}, we comment on the difficulties when extending the ``long resonator'' 
method to the case 
$\sigma=1$. Two key ingredients of the resonance method for a ``long resonator'' are \emph{positivity} and 
\emph{sparsity}. ``Positivity'' means the introduction of an additional function having non-negative Fourier transform 
in $I_1$, which ensures that $I_1$ is a sum of non-negative terms, while ``sparsity'' means that nearby frequencies in 
the resonator function are merged in such a way that one obtains a function having a ``quasi-orthogonality'' property 
which allows to control $I_2$. For details see~\cite{aist,bs1} and the very recent paper~\cite{bs2}. The ``long 
resonator'' has only been implemented for $\sigma$ in the range $[1/2,1)$ so far, since the extension of the method to 
$\sigma=1$ meets serious technical difficulties. The main problem is that the ``sparsification'' of the resonator 
cannot 
be carried out to such an extent as to obtain a truly orthogonal sum, but one rather ends up with a function whose 
square-integral can be estimated only up to multiplicative errors of logarithmic order; if one tried to continue 
thinning out the resonator function to obtain precise control of $I_2$, then from some point on this 
would cause a significant loss in $I_1$ instead. Note that these errors do not play a role in the case $\sigma \in 
[1/2,1)$, where they are negligible in comparison to the main terms, whereas in the case 
$\sigma=1$ we want to obtain a very precise result where only multiplicative errors of order $(1 + 1/\log_2 T)$ are 
allowed.\\

Thus for the case $\sigma=1$ it is necessary to devise a novel variant of the ``long resonator'' technique, which 
will be done in the present paper. The method is genuinely different from those used in~\cite{aist,bs1}, where the 
``sparsification'' of the resonator function played a key role. In the present paper we will avoid this 
sparsity requirement, and actually we make no attempt at all to control the size of $I_2$. Instead we use a 
self-similarity property of the resonator function, which is due to its construction in a completely multiplicative way.

\section{Proof of Theorem~\ref{th1}}

We will use the following approximation formula for the zeta function, which appears in the first lines of the proof 
of Theorem 2 in~\cite{gse}.

\begin{lemma} \label{lemma1}
Define $\zeta(s;Y) := \prod_{p \leq Y} \left(1 - p^{-s} \right)^{-1}$. Let $T$ be large, and set $Y = \exp((\log 
T)^{10})$. Then for $T^{1/10} \leq |t| \leq T$ we have
$$
\zeta(1+it) = \zeta(1+it;Y) \left( 1 + \mathcal{O} \left( \frac{1}{\log T} \right) \right).
$$
\end{lemma}

Instead of this approximation of $\zeta$ by an Euler product we could also use the classical approximation by a 
Dirichlet polynomial (\cite[Theorem 4.11]{titch}), which was used in~\cite{aist,bs1,hilb}. However, the approximation 
by an Euler product is more convenient for us, since the resonator will also be defined as an Euler product.\\

Assume that $T$ is ``large'', and set $Y = \exp((\log T)^{10})$. By Lemma~\ref{lemma1} it suffices to prove Theorem 
\ref{th1} for $\zeta(1+it;Y)$ instead of $\zeta(1+it)$. Set $X = (\log T) (\log_2 T) / 6$, and for primes $p 
\leq X$ set
$$
q_p = \left ( 1 - \frac{p}{X} \right).
$$
This choice of ``weights'' is inspired by those used in the proof of Theorem 2.3 (case $\alpha=1$) in~\cite{hilb}. We 
set $q_1=1$ and $q_p =0$ for $p > X$, and extend our definition in a completely multiplicative way such that we obtain 
weights $q_n$ for all $n \geq 1$. Now we define
$$
R(t) = \prod_{p \leq X} \left(1 - q_p p^{it} \right)^{-1}.
$$
Then 
\begin{eqnarray}
\log (|R(t)|) & \leq & \sum_{p \leq X} \left( \log X - \log p \right) \nonumber\\
& = & \pi(X) \log X - \vartheta(X),
\end{eqnarray}
where $\pi$ is the prime-counting function and $\vartheta$ is the first Chebyshev function. It is well-known that by 
partial summation one has
$$
\pi(X) \log X - \vartheta(X) = \int_2^X \frac{\pi(t)}{t}~ dt = (1 + o(1)) \frac{X}{\log X},
$$
and thus we have
\begin{equation} \label{R_est}
|R(t)|^2 \leq T^{1/3 + o(1)}
\end{equation}
by our choice of $X$.\\

We can write $R(t)$ as a Dirichlet series in the form
\begin{equation} \label{rdi}
R(t) = \sum_{n=1}^\infty q_n n^{it},
\end{equation}
and accordingly
$$
|R(t)|^2 = \sum_{m,n = 1}^\infty q_m q_n \left(\frac{m}{n}\right)^{it}.
$$
Note that all the weights $q_n,~n \geq 1,$ are non-negative reals. Set $\Phi(t)= e^{-t^2}$. By our choice of 
$Y$ we have 
$$
|\zeta(1+it;Y)| \ll \log Y \ll (\log T)^{10}.
$$
Thus, using~\eqref{R_est}, we obtain 
\begin{equation} \label{err1}
\left| \int_{|t| \geq T} \zeta(1+it;Y) |R(t)|^2 \Phi\left(\frac{\log T}{T} t\right) dt \right| \ll 1,
\end{equation}
and 
\begin{equation} \label{err2}
\left| \int_{|t| \leq \sqrt{T}} \zeta(1+it;Y) |R(t)|^2 \Phi\left(\frac{\log T}{T} t\right) dt \right| \ll T^{5/6 + 
o(1)}.
\end{equation}
Using~\eqref{rdi} we can write
\begin{eqnarray*}
I_2 & := & \int_{-\infty}^\infty |R(t)|^2 \Phi\left(\frac{\log T}{T} t\right) dt \\
& = & \sum_{m,n = 1}^\infty 
\int_{-\infty}^{\infty} q_m q_n \left(\frac{m}{n}\right)^{it} \Phi\left(\frac{\log T}{T} 
t\right)~dt.
\end{eqnarray*}
We have the lower bound
\begin{equation} \label{err5}
\int_{\sqrt{T}}^T |R(t)|^2 \Phi\left(\frac{\log T}{T} t\right) dt \gg T^{1+o(1)},
\end{equation}
which follows from $q_1=1$ and the positivity of the Fourier transform of $\Phi$, together with estimates similar to \eqref{err1} and \eqref{err2} to restrict to the desired range. Again using the fact that $\Phi$ has a positive Fourier transform we have
\begin{eqnarray} 
& & \int_{-\infty}^\infty \zeta(1+it;Y) |R(t)|^2 \Phi\left(\frac{\log T}{T} t\right) dt \nonumber\\
& \geq & \int_{-\infty}^\infty \zeta(1+it;X) |R(t)|^2 \Phi\left(\frac{\log T}{T} t\right) dt. \label{xy}
\end{eqnarray}
Here we reduced the range of primes in $\zeta(1+it;~\cdot~)$ from $Y$ to $X$, so that the same primes appear in the 
definitions of $\zeta(1+it;~\cdot~)$ and $R$, respectively. Writing 
$$
\zeta(1+it;X) = \sum_{k=1}^\infty a_k k^{-it}
$$
for appropriate coefficients $(a_k)_{k \geq 1}$, where $a_k \in \{0,1/k\},~k \geq 1$, we have
\begin{eqnarray*}
I_1& := & \int_{-\infty}^\infty \zeta(1+it;X) |R(t)|^2 \Phi\left(\frac{\log T}{T} t\right) dt \\ 
& = & \sum_{k=1}^\infty a_k \sum_{m,n = 1}^\infty
 \int_{-\infty}^{\infty} k^{-it} q_m q_n 
\left(\frac{m}{n}\right)^{it} \Phi\left(\frac{\log T}{T} 
t\right)~dt.
\end{eqnarray*}
Note that we can freely interchange the order of summations and integration, since everything is absolutely 
convergent. Assume $k$ to be fixed. Then, again using the fact that $\Phi$ has a positive Fourier transform, we have
\begin{eqnarray*}
& & \sum_{m,n = 1}^\infty \int_{-\infty}^{\infty} k^{-it} q_m q_n 
\left(\frac{m}{n}\right)^{it} \Phi\left(\frac{\log T}{T} 
t\right) dt \\
& \geq & \sum_{n = 1}^\infty ~\sum_{\substack{1 \leq m < \infty, \\ k | m }} \int_{-\infty}^{\infty} k^{-it} q_m 
q_n 
\left(\frac{m}{n}\right)^{it} \Phi\left(\frac{\log T}{T} 
t\right)dt \\
& = & \sum_{n = 1}^\infty ~\sum_{r=1}^\infty \int_{-\infty}^{\infty} k^{-it} \underbrace{q_{rk}}_{=q_r q_k} 
q_n 
\left(\frac{rk}{n}\right)^{it} \Phi\left(\frac{\log T}{T} 
t\right) dt \\
& = & q_k \underbrace{\sum_{n = 
1}^\infty ~\sum_{r=1}^\infty \int_{-\infty}^{\infty} q_{r} 
q_n 
\left(\frac{r}{n}\right)^{it} \Phi\left(\frac{\log T}{T} 
t\right)~dt.}_{=I_2}
\end{eqnarray*}
Thus we have
\begin{eqnarray}
\frac{I_1}{I_2} & \geq & \sum_{k=1}^\infty a_k q_k \nonumber\\
& = & \prod_{p \leq X} \left(1 - q_p p^{-1}\right)^{-1} \nonumber\\
& = & \left(\prod_{p \leq X} \left(1 - p^{-1} \right)^{-1} \right) \left( \prod_{p \leq X}  \frac{p-1}{p-q_p} 
\right). \label{i1i2}
\end{eqnarray}
For the first product we have
\begin{equation} \label{firstp}
\prod_{p \leq X} \left(1 - p^{-1} \right)^{-1} = e^\gamma \log X + \mathcal{O}(1) = e^\gamma (\log_2 T + 
\log_3 T) + \mathcal{O}(1)
\end{equation}
by Mertens' theorem. For the second product we have
\begin{eqnarray*}
- \log \left(\prod_{p \leq X}  \frac{p-1}{p-q_p} \right) & = & - \sum_{p \leq X} \log \left(1 
- \frac{p}{p+(p-1)X} \right) \\
& \ll & \sum_{p \leq X} \frac{1}{X} \\
& \ll & \frac{1}{\log X}.
\end{eqnarray*}
Thus together with~\eqref{i1i2} and~\eqref{firstp} we have
\begin{eqnarray*}
\frac{I_1}{I_2} \geq e^\gamma (\log_2 T + 
\log_3 T) + \mathcal{O} (1).
\end{eqnarray*}
From~\eqref{err1}--\eqref{xy} we deduce that also
$$
\frac{\left| \int_{\sqrt{T}}^T \zeta(1+it;Y) |R(t)|^2 \Phi\left(\frac{\log T}{T} t\right) dt \right|}{\int_{\sqrt{T}}^T 
|R(t)|^2 
\Phi\left(\frac{\log T}{T} t\right) dt} \geq e^\gamma \left(\log_2 T + \log_3 T \right) + \mathcal{O}(1).
$$
Thus we have
$$
\max_{t \in [\sqrt{T},T]} |\zeta(1+it;Y)| \geq e^\gamma \left(\log_2 T + \log_3 T\right) + \mathcal{O}(1).
$$
As noted at the beginning of the proof, by Lemma~\ref{lemma1} this estimate remains valid if we replace 
$\zeta(1+it;Y)$ by $\zeta(1+it)$. This proves Theorem~\ref{th1}.\\

In conclusion, we add some remarks on the method used in the present paper. As noted in the introduction, earlier 
implementations of the resonance method relied on a combination of ``positivity'' and ``sparsity'' properties. In 
the present paper we show that it is possible to leave out the sparsity requirement, at least in one particular 
instance. One could also obtain the results from~\cite{aist} for the case $\sigma \in (1/2,1)$ using the method from 
the 
present paper without problems. However, it does not seem that the results from~\cite{bs1} for the case $\sigma=1/2$ 
could also be obtained using our argument, since the resonator function there necessarily has a more complicated structure 
(which is not of a simple Euler product form), while our 
argument relies on the fact that one can use a resonator which has completely multiplicative coefficients (which is essential for the argument, since any other coefficients would 
destroy the required self-similarity property).\\

We want to emphasize that the only real restriction on the size of the 
resonator in our argument is~\eqref{R_est}, which gives an upper 
bound for $|R(t)|$. This is different from earlier versions of the ``long resonator'' argument, where bounds on the 
cardinality of the support of $R(t)$ (that is, on the number of non-zero coefficients in the Dirichlet series 
representation) are required.\\

Another remark is that while the sparsity requirement may be unnecessary 
(at least in some cases), the positivity requirement still plays a crucial role in our argument. This prevents a 
possible generalization of the method to the case of functions whose Dirichlet series representation does not only 
contain non-negative real numbers as coefficients. This topic is also discussed in~\cite{MR3554732} in some detail. In 
particular, we have not been able to obtain a result for extreme values of $1/\zeta$ beyond those mentioned in 
\cite{gse}.\\


\end{document}